\documentclass[12pt]{amsart}

\usepackage{DKstyle}
\usepackage{bm}
\usepackage{float}

\newcommand{\Mod}[1]{\ (\mathrm{mod}\ #1)}
\newcommand{\vertiii}[1]{{\left\vert\kern-0.25ex\left\vert\kern-0.25ex\left\vert #1
    \right\vert\kern-0.25ex\right\vert\kern-0.25ex\right\vert}}
 \makeatletter
\newcommand*{\rom}[1]{\expandafter\@slowromancap\romannumeral #1@}
\makeatother
\theoremstyle{plain}

\newtheorem{innercustomthm}{Theorem}
\newenvironment{customthm}[1]
  {\renewcommand\theinnercustomthm{#1}\innercustomthm}
  {\endinnercustomthm}

\usepackage{caption}
\usepackage[labelfont=rm]{subcaption}

\usepackage[pagebackref=true, colorlinks]{hyperref}

\hypersetup{linkcolor=blue,citecolor=blue,urlcolor=blue,menucolor=blue}

\usepackage{comment}


\date{\today}%
\dedicatory{}%
\commby{}%

\title{Quantitative Diophantine approximation with congruence conditions}
\author{Mahbub Alam}

\address{\textbf{Mahbub Alam} \\
    School of Mathematics,
    Tata Institute of Fundamental Research, Mumbai, India 400005 \\
\url{https://sites.google.com/view/mahbubweb}}
\email{mahbub.dta@gmail.com, mahbub@math.tifr.res.in}

\author{Anish Ghosh}

\address{\textbf{Anish Ghosh} \\
    School of Mathematics \\
Tata Institute of Fundamental Research, Mumbai, India 400005}
\email{ghosh@math.tifr.res.in}

\author{Shucheng Yu}
\address{\textbf{Shucheng Yu}\\ Department of Mathematics, Uppsala University, Box 480, SE-75106, Uppsala, SWEDEN}
\email{shucheng.yu@math.uu.se}
\thanks{A.G.\ gratefully acknowledges support from a grant from the Indo-French Centre for the Promotion of Advanced Research, a Department of Science and Technology, Government of India Swarnajayanti fellowship and a MATRICS grant from the Science and Engineering Research Board.
M. A.\ and A. G.\ acknowledge support of the Department of Atomic Energy, Government of India, under project $12-R\&D-TFR-5.01-0500$.
This work received support from a grant from the Infosys foundation.
S.Y.\ acknowledges the support of the Knut and Alice Wallenberg Foundation and ISF grant number  871/17.
S.Y.\ acknowledges that this project has received funding from the European Research Council (ERC) under the European Union's Horizon 2020 research and innovation program (grant agreement No.\ 754475).}
\begin{document}
\begin{abstract}
In this short paper we prove a quantitative version of the Khintchine-Groshev Theorem with congruence conditions. Our argument relies on a classical argument of Schmidt on counting generic lattice points, which in turn relies on a certain variance bound on the space of lattices.
\end{abstract}
\maketitle
\section{Introduction}
Let $\psi : [1, \infty) \to (0, \infty)$ be a continuous and non-increasing function. Let $m, n, d$ be positive integers satisfying $d = m + n$ and let $\vartheta \in \mathrm{M}_{m \times n}(\mathbb{R})$.
Consider the following system of inequalities
\begin{equation}\label{eqdio}
    \norm{\vartheta\bm{q} + \bm{p}}^m < \psi(\norm{\bm{q}}^n),
\end{equation}
with $(\bm{p},\bm{q})\in \Z^m\times \Z^n$. Here $\norm{ \cdot }$ denotes the supremum norm on the corresponding Euclidean spaces. The classical Khintchine-Groshev Theorem gives a criterion on when \eqref{eqdio} has infinitely many integer solutions in $(\bm{p},\bm{q})\in \Z^m\times \Z^n$ for generic (with respect to the Lebesgue measure) $\vartheta\in \textrm{M}_{m\times n}(\R)$:

\begin{customthm}{}[Khintchine-Groshev]\label{thm:kg}
    For almost every (respectively almost no) $\vartheta \in \mathrm{M}_{m \times n}(\mathbb{R})$ there are infinitely many solutions $(\bm{p}, \bm{q}) \in \mathbb{Z}^{m} \times \mathbb{Z}^{n}$ to \eqref{eqdio} if and only if the series $\sum_{t=1}^{\infty}\psi(t)$ diverges (respectively converges).
\end{customthm}


In \cite{NesharimRuhrShi2020} Nesharim-R{\"u}hr-Shi refined the Khintchine-Groshev Theorem by requiring certain congruence conditions: 

\begin{customthm}{}[\textbf{\cite[Theorem 1.2]{NesharimRuhrShi2020}}]
 Let $(\bm{v},N)\in \Z^d\times \N$. Then for almost every $\vartheta \in \mathrm{M}_{m \times n}(\mathbb{R})$ there are infinitely many $(\bm{p}, \bm{q}) \in \mathbb{Z}^{m} \times \mathbb{Z}^{n}$ satisfying $(\bm{p}, \bm{q}) \equiv \bm{v} \Mod{N}$ and \eqref{eqdio} if and only if the series $\sum_{t=1}^{\infty}\psi(t)$ diverges.
\end{customthm}

Schmidt's theorem \cite{Schmidt1960b} on Diophantine approximations refined the Khintchine-Groshev Theorem in another respect, by giving an asymptotic formula for the number of solutions to \eqref{eqdio} as $\|\bm{q}\|$ grows:

\begin{customthm}{}[\textbf{\cite[Theorem 1]{Schmidt1960b}}]
Let $\psi :[1,\infty)\to (0,\infty)$ be as above and let
$$\Psi(T) := \sum_{1\leq t<T} \psi(t) \quad \text{for} \ T \geq 1.$$
Then for almost every $\bm{\vartheta}\in \R^m$, the number of integer solutions in $(\bm{p},q)\in \Z^m\times \Z$ satisfying \eqref{eqdio} (with $n=1$) and $1\leq |q|<T$
is
\begin{align}\label{equ:Schmidtc}
2^d\Psi(T)+O_{\e}\left(\Psi(T)^{1/2+\e}\right).
\end{align}
\end{customthm}
\begin{rmk}
Schmidt actually studied the more refined \textit{one-sided} Diophantine approximation problem by requiring all entries of the vector $(\bm{\vartheta}q+\bm{p}, q)$ to be positive. With this extra requirement one then needs to drop the factor $2^d$ in the above asymptotic formula.
\end{rmk}

In this paper we refine both \cite[Theorem 1.2]{NesharimRuhrShi2020} and \cite[Theorem 1]{Schmidt1960b} by proving a Schmidt-type asymptotic formula for the system of Diophantine inequalities \eqref{eqdio} with extra congruence conditions. 

For the remaining of this paper we fix $d, m,n\in \N$ with $d=m+n$ and we assume $d\geq 3$. We also fix $\psi : [1,\infty)\to (0,\infty)$ a continuous and non-increasing function satisfying that $\sum_{t=1}^{\infty}\psi(t)=\infty$.


We now state our main result.
\begin{Thm}\label{thm:mainresult}
Let $d,m,n\in \N$ and $\psi$ be fixed as above. Let $(\bm{v},N)\in \Z^d\times \N$. 
Let $\nu_1$ and $\nu_2$ be two norms on $\R^m$ and $\R^n$ respectively, with $\nu_2$ normalized such that
\begin{align}\label{equ:norma}
\min_{\bm{w}\in \Z^n\smallsetminus\{\bm{0}\}}\nu_2(\bm{w})=1.
\end{align}
For any $\vartheta\in \mathrm{M}_{m\times n}(\R)$ and for any $T>1$, let $\cN(\vartheta,T)$ 
denote the number of integer solutions in $(\bm{p},\bm{q})\in \Z^m\times \Z^n$ satisfying
\begin{align}\label{equ:dipe}
\nu_1(\vartheta \bm{q}+\bm{p})^m<\psi(\nu_2\left(\bm{q}\right)^n),\ (\bm{p},\bm{q})\equiv \bm{v}\Mod{N},\ 1\leq \nu_2(\bm{q})^n<T.
\end{align}
Then for almost every $\vartheta\in \mathrm{M}_{m\times n}(\R)$,
\begin{align}\label{equ:asyfor}
\cN(\vartheta,T)\ \sim\ N^{-d}c_{\nu_1}c_{\nu_2}\sum_{1\leq t<T}\psi(t) \quad \textrm{as $T\to\infty$},
\end{align}
where for $i=1,2$, $c_{\nu_i}$ is the volume of the unit ball with respect to $\nu_i$.
\end{Thm}
\begin{rmk}\label{rmk:const}
Let $\nu$ be a norm on $\R^{\ell}$ for some $\ell\in \N$.
If $\nu$
is the $L^p$-norm for some $p\geq 1$, then $c_{\nu}=\tfrac{2^{\ell}\Gamma(1+1/p)^{\ell}}{\G(1+\ell/p)}$; 
if $\nu$ is the supremum norm, then $c_{\nu}=2^{\ell}$. In particular, if $N=n=1$ and both $\nu_1$ and $\nu_2$ are the supremum norms, then \thmref{thm:mainresult} recovers the main term in Schmdt's asymptotic formula \eqref{equ:Schmidtc} when $m\geq 2$.
\end{rmk}

Our proof of \thmref{thm:mainresult} consists of two steps: first we prove a counting result regarding the number of points of a generic lattice in an increasing family of Borel sets determined by the inequality system \eqref{equ:dipe}. The second step is then an argument reducing the counting result on generic lattices in the space of lattices to generic lattices in a much smaller sub-manifold of positive co-dimension. We note that this was also the strategy used in \cite{AthreyaParrishTseng2016,AlamGhosh2020,AlamGhosh2020b} where similar quantitative results were proved in various settings for the special approximating function $\psi_0(t)=\frac{c}{t}$ with $c>0$. The arguments used in \cite{AthreyaParrishTseng2016,AlamGhosh2020,AlamGhosh2020b} for the counting result in step one use ergodic theory and are special to $\psi_0$. Here, instead of using this ergodic argument, we use a more soft counting argument of Schmidt \cite{Schmidt1960} which only relies on a variance bound and works for arbitrary increasing family of Borel sets. To incorporate the congruence condition, we use a recently proved variance bound \cite{GhoshKelmerYu2020} on a certain congruence cover of the space of lattices.
\begin{rmk}
We note that the assumption that $d\geq 3$ in \thmref{thm:mainresult} is due to the lack of such variance bounds on congruence covers of the space of rank two lattices. However, when $d=2$ (so that $m=n=1$), $\psi=\psi_0$ and $\nu_1$ and $\nu_2$ are both supremum norms, one can still use the ergodic argument to prove the desired counting result in step one and hence deduce an analogous asymptotic formula in this setting. For simplicity of the presentation, we omit the details here.
\subsection*{Notation and conventions}
\end{rmk}
Throughout the paper, $f(T)\sim g(T)$ means that $\frac{f(T)}{g(T)}\to 1$ as $T\to\infty$. For two positive quantities $A$ and $B$, we will use the notation $A\ll B$ or $A=O(B)$ to mean that there is a constant $c>0$ such that $A\leq cB$, and we will use subscripts to indicate the dependence of the constant on parameters. We will write $A\asymp B$ for $A\ll B\ll A$. For any Borel subset $S$ in a Euclidean space, we use the notation $|S|$ to mean its volume with respect to the usual Lebesgue measure. All vectors in this paper are column vectors even though we will write them as row vectors.

\section{Preliminaries}
\subsection{Siegel transforms and variance estimates}
Let $d\geq 3$ be an integer. Let $G=\SL_d(\R)$ and $\G=\SL_d(\Z)$. It is well known that the homogeneous space $G/\G$ parameterizes $X$, the space of unimodular lattices in $\R^d$ via $g\G\leftrightarrow g\Z^d$. 
More generally, let $(\bm{v},N)\in \Z^d\times \N$ such that $\gcd(\bm{v},N)=1$. Let $X_{N,\bm{v}}$ be the space of affine lattices of the form $g(\Z^n+\tfrac{\bm{v}}{N})$ with $g\in G$. Similar as for $X$, the space $X_{N,\bm{v}}$ can be identified with the homogeneous space $G/\G_{N,\bm{v}}$ via $g\G_{N,\bm{v}} \leftrightarrow g\left(\Z^d+\tfrac{\bm{v}}{N}\right)$, see e.g. \cite[Lemma 3.1]{GhoshKelmerYu2020}. Here
$$\G_{N,\bm{v}}:=\left\{\gamma\in \G: \gamma\bm{v}\equiv \bm{v}\Mod{N}\right\}$$
is the stabilizer in $\G$ of the affine lattice $\Z^d+\tfrac{\bm{v}}{N}$. Note that $\G_{N,\bm{v}}$ is a congruence subgroup (since it contains the principal congruence subgroup $\G(N)$) and $\G_{N,\bm{v}}=\G$ if $N=1$. Generalizing the classical Siegel transform defined on $X$, for any bounded and compactly supported function $f: \R^d\to \C$ we define its \textit{Siegel transform on $X_{N,\bm{v}}\cong G/\G_{N,\bm{v}}$} by
$$\hat{f}\left(g\right):=\sum_{\bm{w}\in \left(\Z^d+\tfrac{\bm{v}}{N}\right)\smallsetminus\{\bm{0}\}}f(g\bm{w}).$$

It was shown in \cite[Proposition 7.1]{MarklofStrombergsson2010} that for any bounded and compactly supported $f$,
\begin{equation}\label{equ:first}
\int_{X_{N,\bm{v}}}\hat{f}(g)d\mu(g)=\int_{\R^d}f(\bm{x})d\bm{x},
\end{equation}
where $\mu$ is the Haar measure of $G$ normalized so that $\mu(X_{N,\bm{v}})=1$. Note that when $f=\chi_{A}$ is the indicator function of some bounded Borel subset $A\subset \R^d$ not containing $\bm{0}$, then
$$\hat{f}(g)=\#(g\left(\Z^d+\tfrac{\bm{v}}{N}\right)\cap A)$$
counts the number of points of the affine lattice $g\left(\Z^d+\tfrac{\bm{v}}{N}\right)$ inside $A$. Hence the integration formula \eqref{equ:first} implies that on average, the counting function $\#(g\left(\Z^d+\tfrac{\bm{v}}{N}\right)\cap A)$ is $|A|$, the volume of $A$. Using \eqref{equ:first} together with a second moment formula \cite[Equation (3.1)]{GhoshKelmerYu2020} the following variance bound was proved in \cite[Corollary 3.4]{GhoshKelmerYu2020}: for any bounded Borel set $A\subset \R^d\smallsetminus \{\bm{0}\}$
\begin{equation}\label{equ:variancebound}
\int_{X_{N,\bm{v}}}\left|\#(g\left(\Z^d+\tfrac{\bm{v}}{N}\right)\cap A)-|A|\right|^2d\mu(g)\ll_{d,N} |A|.
\end{equation}
\subsection{Schmidt's counting results for generic lattices}
It is a classical result by Schmidt \cite[Theorem 1]{Schmidt1960} that given any increasing family of finite-volume Borel sets $\{A_T\}_{T>0}\subset \R^d\smallsetminus\{\bm{0}\}$, (i.e., $A_{T_1}\subset A_{T_2}$ whenever $T_1<T_2$), for $\mu$-a.e.\ unimodular lattice $\Lambda\in X$,
\begin{equation}\label{equ:counting1}
\#(\Lambda\cap A_T)=|A_T|+O_{\e}(|A_T|^{1/2+\e}).
\end{equation}
The main technical tool for Schmidt's arguments is a variance bound in the setting of unimodular lattices, more precisely, the estimate \eqref{equ:variancebound} when $N=1$. In particular, applying Schmidt's arguments and the variance bound \eqref{equ:variancebound} one can get the following counting result for generic lattices in $X_{N,\bm{v}}$, analogous to \eqref{equ:counting1}.

\begin{Prop}\label{prop:schmidtcounting}
Let $(\bm{v},N)\in \Z^d\times \N$ satisfying that $\gcd(\bm{v},N)=1$. Let $\{A_T\}_{T>0}\subset \R^d\smallsetminus\{\bm{0}\}$ be an increasing family of bounded finite-volume Borel subsets. Then for any $\e>0$ and for $\mu$-a.e. $g\in G$, there exists some $T_{g}>0$ such that for all $T>T_{g}$
\begin{align*}
\#\left(g(\Z^d+\tfrac{\bm{v}}{N})\cap A_T\right)=|A_T|+O_{\e}\left(|A_T|^{1/2+\e}\right).
\end{align*}
\end{Prop}
Using a simple scaling argument we have the following counting result which we will use later.
\begin{Cor}\label{cor:schmidtcounting}
Let $(\bm{v},N)\in \Z^d\times \N$ and let $\{A_T\}_{T>0}$ be as in \propref{prop:schmidtcounting}. We further assume that $\lim\limits_{T\to\infty}|A_T|=\infty$. Then for $\mu$-a.e. $g\in G$,
$$\#(g(N\Z^d+\bm{v})\cap A_T)\ \sim\ N^{-d}|A_T|\quad \textrm{as $T\to\infty$}.$$
\end{Cor}
\begin{proof}
Let $l=\gcd(\bm{v},N)$ and let $\bm{v}'=\bm{v}/l$ and $N'=N/l$ so that $\tfrac{\bm{v}}{N}=\tfrac{\bm{v}'}{N'}$ and $\gcd(\bm{v}',N')=1$. Applying \propref{prop:schmidtcounting} to the space $X_{N',\bm{v}'}$ and the family $\left\{N^{-1}A_T\right\}_{T>0}$, (here we can apply \propref{prop:schmidtcounting} since $\lim\limits_{T\to\infty}|N^{-1}A_T|=\lim\limits_{T\to\infty}N^{-d}|A_T|=\infty$) to get for $\mu$-a.e. $g\in G$,
\begin{align*}
\#\left(g(\Z^d+\tfrac{\bm{v}'}{N'})\cap N^{-1}A_T\right)\ \sim\ |N^{-1}A_T|=N^{-d}|A_T|\quad\textrm{as $T\to\infty$}.
\end{align*}
We can thus finish the proof by noting that
\begin{displaymath}
\#\left(g(\Z^d+\tfrac{\bm{v}'}{N'})\cap N^{-1}A_T\right)=\#\left(g(\Z^d+\tfrac{\bm{v}}{N})\cap N^{-1}A_T\right)=\#\left(g\left(N\Z^d+\bm{v}\right)\cap A_T\right).\qedhere
\end{displaymath}
\end{proof}
\subsection{Relating to counting lattice points}
Let $U := u(\mathrm{M}_{m \times n}(\mathbb{R}))< G$, where $u : \mathrm{M}_{m \times n}(\mathbb{R}) \to G$ is defined as
\[
    u(\vartheta) := \begin{pmatrix}
        1_m & \vartheta \\
        0   & 1_n
    \end{pmatrix}.
\]
Here $1_\ell$ denotes the $\ell \times \ell$ identity matrix and $0$ is the zero matrix in $\textrm{M}_{n\times m}(\R)$. Note that
$$N\Z^d+\bm{v}=\left\{(\bm{p},\bm{q})\in \Z^m\times \Z^n: (\bm{p},\bm{q})\equiv \bm{v}\Mod{N}\right\}$$
and any lattice point of $u(\vartheta)\left(N\Z^d+\bm{v}\right)$ is of the form
$$u(\vartheta)(\bm{p},\bm{q})=(\vartheta \bm{q}+\bm{p}, \bm{q})$$
for some $(\bm{p},\bm{q})\in N\Z^d+\bm{v}$. This, together with the normalization \eqref{equ:norma} 
implies that
\begin{align}\label{equ:relation}
\cN(\vartheta,T)=\#\left(u(\vartheta)\left(N\Z^d+\bm{v}\right)\cap E_T)\right),
\end{align}
where for any $T>1$
\begin{align}\label{equ:set}
E_T=E_{\psi,\nu_1,\nu_2,T}:=\left\{(\bm{x},\bm{y})\in \R^m\times\R^n: \nu_1(\bm{x})^m<\psi\left(\nu_2(\bm{y})^n\right),\ 1\leq \nu_2(\bm{y})^n<T\right\}.
\end{align}
\subsection{Volume calculation}\label{sec:voulme}
In this subsection we give a quick computation for the volume of the sets $E_T$ defined in \eqref{equ:set}. We first record a simple volume formula for balls with respect to a norm $\nu$ on $\R^{\ell}$. For all $r>0$ let $B_{\nu}(r):=\{\bm{x}\in \R^{\ell}: \nu(\bm{x})<r\}$ be the open radius $r$-ball with respect to $\nu$, centered at the origin. Since $\nu$ is \textit{positive homogeneous}, i.e., $\nu(r\bm{x})=r\nu(\bm{x})$ for any $r>0$ and $\bm{x}\in \R^{\ell}$, we have
\begin{align}\label{equ:ball}
|B_{\nu}(r)|=|rB_{\nu}(1)|=|B_{\nu}(1)|r^{\ell}=c_{\nu}r^{\ell}.
\end{align}
Now using the identity \eqref{equ:ball} we have
\begin{align*}
|E_T|&=\int_{\{\bm{y}\in \R^n: 1\leq \nu_2(\bm{y})^n<T\}}\int_{\{\bm{x}\in \R^m: \nu_1(\bm{x})^m<\psi(\nu_2(\bm{y})^n)\}}d\bm{x}d\bm{y}\\
&=c_{\nu_1}\int_{\{\bm{y}\in \R^n: 1\leq \nu_2(\bm{y})^n<T\}}\psi(\nu_2(\bm{y})^n)d\bm{y}=:c_{\nu_1}F(T).
\end{align*}
It is not hard to check that (again using \eqref{equ:ball}) for any $T>1$, $F'(T)=c_{\nu_2}\psi(T)$. This, together with the fact that $F(1)=0$, implies that $F(T)=c_{\nu_2}\int_{1}^T\psi(r)dr$. Hence
\begin{align}\label{equ:volume}
|E_T|=c_{\nu_1}c_{\nu_2}\int_{1}^T\psi(r)dr=c_{\nu_1}c_{\nu_2}\sum_{1\leq t<T}\psi(t)+O_{\psi}(1).
\end{align}

\subsection{Decomposition of the Haar measure}
Let $H< G$ be the parabolic subgroup such that
\begin{align*}
H:=\left\{h=\begin{pmatrix}
\alpha & 0\\
\beta & \gamma \end{pmatrix}\in G: \alpha\in \GL_m(\R),\ \gamma\in \GL_n(\R),\ \beta\in \textrm{M}_{n\times m}(\R)\right\}.
\end{align*}
Here $0$ denotes the zero matrix in $\textrm{M}_{m\times n}(\R)$. We note that there is a Zariski dense subset of $G$ such that any $g$ in this subset can be written uniquely as the product $g=hu(\vartheta)$ with $h\in H$ and $u(\vartheta)\in U$. We note that under this decomposition, the Haar measure $\mu$ decomposes as (up to scalars) $d\mu(g)=dhd\vartheta$, where $dh$ is a left $H$-invariant Haar measure of $H$ and $d\vartheta$ is the usual Lebesgue measure on $U(\cong \textrm{M}_{m\times n}(\R)\cong \R^{mn})$. In view of this measure decomposition and Fubini's theorem we can restate \corref{cor:schmidtcounting}:
\begin{Prop}\label{prop:countingfinal}
Keep the assumptions as in \corref{cor:schmidtcounting}. Then for a.e. $h\in H$ (with respect to $dh$) and for a.e. $\vartheta\in \textrm{M}_{m\times n}(\R)$ (with respect to $d\vartheta$)
$$\#(hu(\vartheta)(N\Z^d+\bm{v})\cap A_T)\ \sim\ N^{-d}|A_T|\quad \textrm{as $T\to\infty$}.$$
\end{Prop}
\section{Proof of \thmref{thm:mainresult}}
In this section we give the proof of \thmref{thm:mainresult}. The following lemma is the key step for our reduction argument. It says that the sets $E_T$ are stable under small perturbations of elements in $H$ close to the identity matrix.

\begin{Lem}\label{lem:cont}
There exists $0<c_0<\frac{1}{2}$ such that for any $\e\in (0,c_0)$ there exists an open neighborhood $H_{\e}\subset H$ of the identity element such that for all $h\in H_{\e}$
$$E_{T,\e}^{-}\subset hE_T\subset E_{T,\e}^+,\ \forall\ T>10,$$
where $E_T$ is as given in \eqref{equ:set} and
$$E_{T,\e}^{-}:=\left\{(\bm{x},\bm{y})\in \R^m\times \R^n: \nu_1(\bm{x})^m<(1+\e)^{-1}\psi\left((1+\e)\nu_2(\bm{y})^n\right), \frac{3}{2}\leq \nu_2(\bm{y})^n<(1+\e)^{-1}T \right\},$$
and $E_{T,\e}^+:=E'_{T,\e}\cup C_0$ with
$$E_{T,\e}':=\left\{(\bm{x},\bm{y})\in \R^m\times \R^n: \nu_1(\bm{x})^m<(1+\e)\psi\left((1+\e)^{-1}\nu_2(\bm{y})^n\right), \frac{3}{2}\leq \nu_2(\bm{y})^n<(1+\e)T \right\},$$
and
$$C_0:=\left\{(\bm{x},\bm{y})\in \R^m\times \R^n: \nu_1(\bm{x})^m<2\psi(1),\ 
\frac{1}{2}<\nu_2(\bm{y})^n\leq \frac{3}{2}\right\}.$$
\end{Lem}
\begin{proof}
View $(\R^m,\nu_1)$ and $(\R^n,\nu_2)$ as two normed spaces, and for any $\alpha\in \GL_m(\R)$, $\beta\in \textrm{M}_{m\times n}(\R)$ and $\gamma\in \GL_n(\R)$, we denote by $\|\alpha\|_{\nu_1}$, $\|\beta\|_{\nu_1,\nu_2}$ and $\|\gamma\|_{\nu_2}$ their corresponding operator norms. For each $\e>0$ let
\begin{align*}
\tilde{H}_{\e}:=\left\{h=\begin{pmatrix}
\alpha & 0\\
\beta & \gamma \end{pmatrix}\in G: \max\{\|\alpha\|^m_{\nu_1},\|\gamma\|^n_{\nu_2}\}<1+\frac{\e}{2},\ \|\beta\|_{\nu_1,\nu_2}<\frac{\e}{4n\psi(1)^{1/m}} \right\}.
\end{align*}
and define $H_{\e}:=\tilde{H}_{\e}\cap \tilde{H}_{\e}^{-1}$. Then clearly $H_{\e}\subset H$ is an open neighborhood of the identity element. We need to prove the above inclusion relations for all $h\in H_{\e}$.

Fix $h=\left(\begin{smallmatrix}
\alpha & 0\\
\beta & \gamma\end{smallmatrix}\right)\in H_{\e}$ and $T>10$. We first prove the relation $hE_T\subset E_{T, \epsilon}^+ = E_{T,\e}'\cup C_0$. Let $(\bm{x},\bm{y})\in E_T$. First note that $\nu_1(\bm{x})^m<\psi(\nu_2(\bm{y})^n)\leq \psi(1)$, where for the second inequality we used the assumption that $\nu_2(\bm{y})\geq 1$ and that $\psi$ is non-increasing. Now using the definition of $\tilde{H}_{\e}$ and the fact that $h\in H_{\e}\subset \tilde{H}_{\e}$ we have
$$h(\bm{x},\bm{y})=(\alpha\bm{x}, \beta\bm{x}+\gamma\bm{y})$$
with
\begin{align*}
\nu_1(\alpha\bm{x})^m\leq \left(1+\tfrac{\e}{2})\nu_1(\bm{x}\right)^m
\end{align*}
and (using the triangle inequality and the inequalities that $\nu_2(\bm{x})^m<\psi(1)$ and $\nu_2(\bm{y})>1$)
\begin{align*}
\nu_2(\beta\bm{x}+\gamma\bm{y})\leq \nu_2(\gamma\bm{y})+\nu_2(\beta\bm{x})\leq (1+\tfrac{\e}{2})^{1/n}\nu_2(\bm{y})+\frac{\e}{4n\psi(1)^{1/m}}\nu_1(\bm{x})<(1+\e)^{1/n}\nu_2(\bm{y}).
\end{align*}
Here for the last inequality we used the inequality that $(1+\tfrac{\e}{2})^{1/n}+\frac{\e}{4n}<(1+\e)^{1/n}$ which can be guaranteed for all $\epsilon\in (0,c_0)$ by taking $c_0>0$ sufficiently small.

Now if $\nu_2(\beta\bm{x}+\gamma\bm{y})^n>\frac{3}{2}$, then
$$\nu_2(\beta\bm{x}+\gamma\bm{y})^n<(1+\e)\nu_2(\bm{y})^n<(1+\e)T,$$
and
\begin{align*}
\nu_1(\alpha\bm{x})^m\leq (1+\tfrac{\e}{2})\nu_1(\bm{x})^m<(1+\tfrac{\e}{2})\psi(\nu_2(\bm{y})^n)<(1+\e)\psi\left(\left(1+\e)^{-1}\nu_2(\beta\bm{x}+\gamma\bm{y})^n\right)\right).
\end{align*}
This implies that $h(\bm{x},\bm{y})\in E'_{T,\e}\subset E_{T,\e}^{+}$. If $\nu_2(\beta\bm{x}+\gamma\bm{y})\leq \frac{3}{2}$, then
\begin{align*}
\nu_1(\alpha\bm{x})^m\leq (1+\tfrac{\e}{2})\nu_1(\bm{x})^m<2\psi(1),
\end{align*}
and (since $h^{-1}=\left(\begin{smallmatrix}
\alpha^{-1} & 0\\
-\gamma^{-1}\beta\alpha^{-1} & \gamma^{-1}\end{smallmatrix}\right)\in \tilde{H}_{\e}$)
\begin{align*}
\nu_2(\beta\bm{x}+\gamma\bm{y})\geq \nu_2(\gamma\bm{y})-\nu_2(\beta\bm{x})\geq (1+\tfrac{\e}{2})^{-1/n}\nu_2(\bm{y})-\frac{\e}{4n}\geq (1+\tfrac{\e}{2})^{-1/n}-\frac{\e}{4n}>\frac{1}{2^{1/n}},
\end{align*}
where the last inequality again can be guaranteed by taking $c_0$ sufficiently small. Thus in this case we have $h(\bm{x},\bm{y})\in C_0\subset E_{T,\e}^+$. This finishes the proof of the relation $hE_T\subset E_{T,\e}^+$. Similarly, using the fact that $h^{-1},h\in \tilde{H}_{\e}$, one can show $h^{-1}E_{T,\e}^{-}\subset E_T$, or equivalently, $E_{T,\e}^{-}\subset hE_T$.
\end{proof}
\begin{rmk}\label{rmk:volume}
Using similar volume calculations as in \secref{sec:voulme} one can easily see that for all $T>10$ and $\e\in (0,c_0)$
\begin{align}\label{equ:volume2}
|E_{T,\e}^{\pm}|&=(1+ \e)^{\pm 1}c_{\nu_1}c_{\nu_2}\int_{\frac{3}{2}}^{(1+\e)^{\pm 1}T}\psi((1+\e)^{\mp 1}r)dr+O(1)\nonumber\\
&=(1+ \e)^{\pm 2}c_{\nu_1}c_{\nu_2}\int_{1}^{T}\psi(r)dr+O_{\psi,\nu_1,\nu_2}(1).
\end{align}
\end{rmk}
\begin{proof}[Proof of \thmref{thm:mainresult}]
For simplicity of notation, let us denote $\Lambda_0:=N\Z^d+\bm{v}$. Fix a sequence of positive numbers $\{\e_{\ell}\}_{\ell\in\N}\subset (0,c_0)$ with $\e_{\ell}\to 0$ as $\ell\to\infty$. Here $c_0$ is the constant as in \lemref{lem:cont}. Note that for each $\ell\in \N$, both $\{E_{T,\e_{\ell}}^{\pm}\}_{T>10}$ are increasing families of bounded Borel subsets not containing $\bm{0}$. Moreover, using the volume formulas \eqref{equ:volume} and \eqref{equ:volume2} 
we have
\begin{align}\label{equ:volest}
\frac{|E_{T,\e_{\ell}}^+|}{|E_T|}\leq (1+\e_{\ell})^2+O_{\psi,\nu_1,\nu_2}(|E_T|^{-1})\quad \textrm{and}\quad \frac{|E_{T,\e_{\ell}}^-|}{|E_T|}\geq (1+\e_{\ell})^{-2}+O_{\psi,\nu_1,\nu_2}(|E_T|^{-1}).
\end{align}
In particular, the above second estimate, together with the relation $E_T\subset E_{T,\e_{\ell}}^+$ and the fact that
$$|E_T|\asymp \sum_{1\leq t<T}\psi(t)\to \infty\quad \textrm{as $T\to\infty$},$$
implies that $\lim\limits_{T\to\infty}|E_{T,\e_{\ell}}^{\pm}|=\infty$. Thus we can apply \propref{prop:countingfinal} for the two families $\{E_{T,\e_{\ell}}^{\pm}\}_{T>10}$. In fact, combining \propref{prop:countingfinal} and \lemref{lem:cont} and using the fact that a finite intersection of full measure sets is still of full measure (hence intersecting any open set nontrivially) we can find, for each $\ell\in \N$, $h_{\ell}\in H$ such that
\begin{align}\label{equ:re1}
E_{T,\e_{\ell}}^{-}\subset h_{\ell}E_T\subset E_{T,\e_{\ell}}^+,\ \forall\ T>10,
\end{align}
and that there exists a full measure subset $U_{\ell}\subset U$ such that for all $u(\vartheta)\in U_{\ell}$
\begin{align}\label{equ:re2}
\#(h_{\ell}u(\vartheta)\Lambda_0\cap E_{T,\e_{\ell}}^{\pm})\ \sim\ N^{-d}|E_{T,\e_{\ell}}^{\pm}|\quad \textrm{as $T\to\infty$}.
\end{align}
Now let $U_{\infty}=\cap_{\ell\in \N}U_{\ell}$ which is still of full measure. In view of the relation \eqref{equ:relation} and the volume calculation \eqref{equ:volume}, 
it suffices to show that for all $\vartheta\in U_{\infty}$
\begin{align*}
\#\left(u(\vartheta)\Lambda_0\cap E_T)\right)\ \sim\ N^{-d}|E_T|\quad \textrm{as $T\to\infty$}.
\end{align*}
Now for any $u(\vartheta)\in U_{\infty}$ using the relation \eqref{equ:re1} we have for each $\ell\in \N$
\begin{align*}
\#(h_{\ell}u(\vartheta)\Lambda_0\cap E_{T,\e_{\ell}}^{-})\leq \#(h_{\ell}u(\vartheta)\Lambda_0\cap h_{\ell} E_T)\leq \#(h_{\ell}u(\vartheta)\Lambda_0\cap E_{T,\e_{\ell}}^+),
\end{align*}
or equivalently,
\begin{align*}
\frac{\#(h_{\ell}u(\vartheta)\Lambda_0\cap E_{T,\e_{\ell}}^{-})}{N^{-d}|E_T|}\leq \frac{\#(u(\vartheta)\Lambda_0\cap E_T)}{N^{-d}|E_T|}\leq \frac{\#(h_{\ell}u(\vartheta)\Lambda_0\cap E_{T,\e_{\ell}}^+)}{N^{-d}|E_T|}.
\end{align*}
Now applying \eqref{equ:re2} and \eqref{equ:volest} we get
\begin{align*}
\limsup_{T\to\infty}\frac{\#(u(\vartheta)\Lambda_0\cap E_T)}{N^{-d}|E_T|}\leq \liminf_{T\to\infty}\frac{|E_{T,\e_{\ell}}^+|}{|E_T|}\leq (1+\e_{\ell})^2,
\end{align*}
and
\begin{align*}
\liminf_{T\to\infty}\frac{\#(u(\vartheta)\Lambda_0\cap E_T)}{N^{-d}|E_T|}\geq \limsup_{T\to\infty}\frac{|E_{T,\e_{\ell}}^-|}{|E_T|}\geq (1+\e_{\ell})^{-2}.
\end{align*}
Finally taking $\ell\to\infty$ finishes the proof.
\end{proof}
\bibliographystyle{alpha}
\bibliography{DKbibliog}

\begin{thebibliography}{GKY20}

\bibitem[AG20a]{AlamGhosh2020}
M.~Alam and A.~Ghosh.
\newblock Equidistribution on homogeneous spaces and the distribution of
  approximates in {D}iophantine approximation.
\newblock {\em Trans. Amer. Math. Soc.}, 373(5):3357--3374, 2020.

\bibitem[AG20b]{AlamGhosh2020b}
M.~Alam and A.~Ghosh.
\newblock Quantitative rational approximation on spheres.
\newblock {\em arXiv preprint arXiv:2003.02243}, 2020.

\bibitem[APT16]{AthreyaParrishTseng2016}
J.~Athreya, A.~Parrish, and J.~Tseng.
\newblock Ergodic theory and {D}iophantine approximation for translation
  surfaces and linear forms.
\newblock {\em Nonlinearity}, 29(8):2173--2190, 2016.

\bibitem[GKY20]{GhoshKelmerYu2020}
A.~Ghosh, D.~Kelmer, and S.~Yu.
\newblock {Effective Density for Inhomogeneous Quadratic Forms I: Generic Forms
  and Fixed Shifts}.
\newblock {\em Int. Math. Res. Not. IMRN}, 08 2020.
\newblock rnaa206.

\bibitem[MS10]{MarklofStrombergsson2010}
J.~Marklof and A.~Str\"{o}mbergsson.
\newblock The distribution of free path lengths in the periodic {L}orentz gas
  and related lattice point problems.
\newblock {\em Ann. of Math. (2)}, 172(3):1949--2033, 2010.

\bibitem[NRS20]{NesharimRuhrShi2020}
E.~Nesharim, R.~R\"{u}hr, and R.~Shi.
\newblock Metric {D}iophantine approximation with congruence conditions.
\newblock {\em Int. J. Number Theory}, 16(9):1923--1933, 2020.

\bibitem[Sch60a]{Schmidt1960b}
M.~W. Schmidt.
\newblock A metrical theorem in diophantine approximation.
\newblock {\em Canadian J. Math.}, 12:619--631, 1960.

\bibitem[Sch60b]{Schmidt1960}
W.~M. Schmidt.
\newblock A metrical theorem in geometry of numbers.
\newblock {\em Trans. Amer. Math. Soc.}, 95:516--529, 1960.

\end{thebibliography}

\end{document}